\documentclass{article}
\renewcommand{\baselinestretch}{1.18}
\usepackage{amssymb}

\newcommand{\eqref}[1]{\mbox{\rm(\ref{#1})}}
\newcommand{\claim}[1]{\smallskip\noindent{\bf Claim #1} }
\newcommand{\rcr}[1]{\ref{cr:\pl:#1}}

\newcommand{\imply}{\Rightarrow}
\newcommand{\pl}{}
\newcommand{\C}[1]{{\protect\cal #1}}
\newcommand{\B}[1]{{\bf #1}}
\newcommand{\I}[1]{{\mathbb #1}}
\renewcommand{\O}[1]{\overline{#1}}
\newcommand{\binom}[2]{{#1\choose #2}}
\newcommand{\e}{\epsilon}
\newcommand{\comment}[1]{}
\newcommand{\qed}{\nolinebreak\mbox{\hspace{5 true pt}%
\rule[-0.85 true pt]{3.9 true pt}{8.1 true pt}}}

\newtheorem{theorem}{Theorem}
\newtheorem{lemma}[theorem]{Lemma}
\newtheorem{corollary}[theorem]{Corollary}

\newtheorem{problem}[theorem]{Problem}

\begin{document}

\newcommand{\Rbinom}[2]{\binom{#1}{#2}}
\newcommand{\R}{\Pi}
\newcommand{\sol}{sol}
\newcommand{\dmin}{d_0} %{\mbox{\rm\tiny m}}}
\newcommand{\dminmin}{d'} %{\mbox{\rm\tiny m}}}
\newcommand{\tmin}{t_0} %{\mbox{\rm\tiny m}}}
\newcommand{\vmax}{v_0} %{\mbox{\rm\tiny m}}}
\newcommand{\fmax}{f_0} %{\mbox{\rm\tiny m}}}
\newcommand{\need}[1]{\par\noindent{\bf We need here: $#1$.}\par}
\newcommand{\cone}{m_1}
\newcommand{\ctwo}{m_2}
\newcommand{\cthree}{m_3}

\newcommand{\bpr}[1]{
 \renewcommand{\thetheorem}{#1}
 \begin{problem}\label{pr:\pl:#1}}

\newcommand{\epr}{
 \end{problem}
 \renewcommand{\thetheorem}{\arabic{theorem}}
} 

\newcommand{\rpr}[1]{\ref{pr:\pl:#1}}

\title{Asymptotic Size Ramsey Results for Bipartite Graphs}
\author{Oleg Pikhurko\thanks{Supported by a Research Fellowship, St.\
John's College, Cambridge. Part of this research was carried out
during the author's stay at the Humboldt University, Berlin,
sponsored by the German Academic Exchange Service (DAAD).}\\
 DPMMS, Centre for Mathematical Sciences\\
 Cambridge University, Cambridge~CB3~0WB, England\\
 E-mail: {\tt O.Pikhurko@dpmms.cam.ac.uk}}
\maketitle

\begin{abstract} We show that $\lim_{n\to\infty}\hat
r(F_{1,n},\dots,F_{q,n},F_{p+1},\dots,F_{r})/n$ exists, where the
bipartite graphs $F_{q+1},\dots,F_r$ do not depend on $n$ while, for
$1\le i\le q$, $F_{i,n}$ is obtained from some bipartite graph
$F_i$ with parts $V_1\cup V_2=V(F_i)$ by duplicating each vertex $v\in
V_2$ $(c_v+o(1))n$ times for some real $c_v>0$.

In fact, the limit is the minimum of a certain mixed integer
program. Using the Farkas Lemma we compute it when each forbidden
graph is a complete bipartite graph, in particular answering a
question of Erd\H os, Faudree, Rousseau and Schelp~(1978) who asked
for the asymptotics of $\hat r(K_{s,n},K_{s,n})$ for fixed $s$ and
large $n$. Furthermore, we prove (for all sufficiently large $n$) the
conjecture of Faudree, Rousseau and Sheehan~(1983) that $\hat
r(K_{2,n},K_{2,n}) =18n-15$.\end{abstract}

\section{Introduction}

Let $(F_1,\dots,F_r)$ be an $r$-tuple of graphs which are called {\em
forbidden}.  We say that a graph $G$ {\em arrows} $(F_1,\dots,F_r)$ if
for any $r$-colouring of $E(G)$, the edge set of $G$, there is a copy
of $F_i$ of colour $i$ for some $i\in[r]:=\{1,\dots,r\}$. We denote
this {\em arrowing property} by $G\to(F_1,\dots,F_r)$.

The (ordinary) {\em Ramsey number} asks for the minimum order of such
$G$. Here, however, we deal exclusively with the {\em size Ramsey
number}
 $$\hat r(F_1,\dots,F_r)=\min\{e(G)\mid G\to(F_1,\dots,F_r)\}$$
 which is the smallest number of edges that an arrowing graph
can have.

Size Ramsey numbers seem hard to compute, even for simple forbidden
graphs. For example, the old conjecture of Erd\H os~\cite{erdos:81}
that $\hat r(K_{1,n},K_3)=3n(n+1)/2$ has only recently been disproved
in~\cite{pikhurko:99c}, where it is shown that $\hat
r(K_{1,n},F)=(1+o(1))n^2$ for any fixed $3$-chromatic graph
$F$. (Here, $K_{m,n}$ is the complete bipartite graph with parts of
sizes $m$ and $n$; $K_n$ is the complete graph of order $n$.)

This research initiated as an attempt to find the asymptotics of $\hat
r(K_{1,n}, F)$ for a fixed graph $F$. The case $\chi(F)\ge 4$ is
treated in~\cite{pikhurko:01} (and~\cite{pikhurko:99c} deals with
$\chi(F)=3$). What can be said if $F$ is a bipartite graph?

Faudree, Rousseau and Sheehan~\cite{faudree+rousseau+sheehan:83}
proved that
 $$\hat r(K_{1,n},K_{2,m})=4n+2m-4$$
 for every $m\ge9$ if $n$ is sufficiently large (depending on $m$) and
stated that their method shows that $\hat r(K_{1,n},K_{2,2})=4n$.
They also observed that $K_{s,2n}$ arrows the pair $(K_{1,n},C_{2s})$,
where $C_{2s}$ is a cycle of order $2s$; hence $\hat
r(K_{1,n},C_{2s})\le 2sn$.

Let $P_s$ be the path with $s$ vertices. Lortz and
Mengersen~\cite{lortz+mengersen:98} showed that $K_{k,2n-1}\to
(K_{1,n},P_{2k+1})$ and $K_k+\O{K}_{2n-k-1}\to (K_{1,n},P_{2k})$ and
conjectured that this is sharp for any $s\ge4$ provided $n$ is
sufficiently large, that is,
 \begin{equation}\label{eq:\pl:lm97}\hat r(K_{1,n},P_s)=\left\{\begin{array}{ll}
 2kn-k,&\mbox{if $s=2k+1$,}\\
 2kn-k(k+3)/2,&\mbox{if $s=2k$,}\end{array}\right.\quad\quad  n\ge n_0(s).\end{equation}
 The conjecture was proved for $4\le s\le 7$ in~\cite{lortz+mengersen:98}.

Size Ramsey numbers $\hat r(F_1,F_2)$ for bipartite graphs $F_1$ and
$F_2$ (and in some papers $F_1$ is a small star) are also studied
in~\cite{erdos+faudree+rousseau+schelp:78,burr+erdos+faudree+rousseau+schelp:78,beck:83,beck:90,erdos+faudree:92,erdos+rousseau:93,ke:93,haxell+kohayakawa:95}
for example.

It is not hard to see that, for fixed $s_1,\dots,s_r\in \I N$ and
$t_1,\dots,t_r\in\I R_{>0}$, we have
 \begin{equation}\label{eq:\pl:On}\hat r(K_{s_1,\lfloor t_1n\rfloor},\dots,K_{s_r,\lfloor
t_rn\rfloor}) = O(n).\end{equation}
 This follows, for example, by assuming that $s_1=\dots=s_r=s$,
$t_1=\dots=t_r=t$ and considering $K_{v_1,v_2}$, where $v_1=(s-1)r+1$
and $v_2=\lceil rtn\binom {v_1}s\rceil$. The latter graph has the
required arrowing property. Indeed, for any $r$-colouring, each vertex
of $V_2$ is incident to at least $s$ edges of same colour; hence there
are at least $v_2$ monochromatic $K_{s,1}$-subgraphs and some $S\in
\binom{V_1}s$ appears in at least $rtn$ such subgraphs of which at
least $tn$ have same colour.

Here we will show that the limit $\lim_{n\to\infty} \hat r
(F_{1,n},\dots,F_{r,n})/n$ exists if each forbidden graph is either a
fixed bipartite graph or a subgraph of $K_{s,\lfloor tn\rfloor}$ which
`dilates' uniformly with $n$ (the precise definition will be given in
Section~\ref{dilatation}). In particular, $\hat r(K_{1,n},F)/n$ tends
to a limit for any fixed bipartite graph $F$.

The limit value can in fact be obtained as the minimum of a certain
mixed integer program (which does depend on $n$).  We have been able
to solve the MIP when each $F_{i,n}$ is a complete bipartite
graph. In particular, this answers a question by Erd\H os, Faudree,
Rousseau and Schelp~\cite[Problem~B]{erdos+faudree+rousseau+schelp:78}
who asked for the asymptotics of $\hat r(K_{s,n},K_{s,n})$.  Working
harder on the case $s=2$ we prove (for all sufficiently large $n$) the
conjecture of Faudree, Rousseau and
Sheehan~\cite[Conjecture~15]{faudree+rousseau+sheehan:83} that
 \begin{equation}\label{eq:\pl:frs83} \hat r(K_{2,n},K_{2,n})=18n-15,\end{equation}
 where the upper bound is obtained by considering
$K_{3,6n-5}\to(K_{2,n},K_{2,n})$. Unfortunately, the range on $n$
from~\eqref{eq:\pl:frs83} is not specified
in~\cite{faudree+rousseau+sheehan:83}, although it is stated there
that $\hat r(K_{2,2},K_{2,2})=15$, where the upper bound follows
apparently from $K_6\to (K_{2,2},K_{2,2})$.

Unfortunately, our MIP is not well suited for practical calculations
and we were not able to compute the asymptotics for any other
non-trivial forbidden graphs; in particular, we had no progress
on~\eqref{eq:\pl:lm97}. But we hope that the introduced method will produce
more results: although the MIP is hard to solve, it may well be
possible that, for example, some manageable relaxation of it gives
good lower or upper bounds.

Our method does not work if we allow both vertex classes of forbidden
graphs to grow with $n$. In these settings, in fact, we do not know
the asymptotics even in simplest cases. For example, the best known
bounds on $r=\hat r(K_{n,n},K_{n,n})$ seem to be $r<\frac32n^32^n$ for
$n\ge 6$ (Erd\H os, Faudree, Rousseau and
Schelp~\cite{erdos+faudree+rousseau+schelp:78}) and $r>
\frac1{60}n^22^n$ (Erd\H os and Rousseau~\cite{erdos+rousseau:93}).

Our theorem on the existence of the limit can be extended to generalized
size Ramsey problems; this is discussed in Section~\ref{generalized}.

\section{Some Definitions}\label{dilatation}

We decided to gather most of the definitions in this section for quick
reference.

We assume that bipartite graphs come equipped with a fixed bipartition
$V(F)=V_1(F)\cup V_2(F)$, although graph embeddings need not preserve
it. We denote $v_i(F)=|V_i(F)|$, $i=1,2$; thus $v(F)=v_1(F)+v_2(F)$.

For $A\subset V_1(F)$, we define 
 $$F^A=\{v\in V_2(F)\mid \Gamma_F(v)=A\},$$
 where $\Gamma_F(v)$ denotes the neighbourhood of $v$ in $F$. (We will
write $\Gamma(v)$, etc., when the encompassing graph $F$ is clear from
the context.) Clearly, in order to determine $F$ (up to an
isomorphism) it is enough to have $V_1(F)$ and $|F^A|$ for all
$A\subset V_1(F)$. This motivates the following definitions.

A {\em weight} $\B f$ on a set $V(\B f)$ a sequence $(f_A)_{A\in
2^{V(\B f)}}$ of non-negative reals. A bipartite graph $F$ {\em
agrees} with $\B f$ if $V_1(F)=V(\B f)$ and $F^A=\emptyset$ if and
only if $f_A=0$, $A\in 2^{V(\B f)}$. A sequence of bipartite graphs
$(F_n)_{n\in\I N}$ is a {\em dilatation} of $\B f$ (or {\em dilates}
$\B f$) if each $F_n$ agrees with $\B f$ and
 $$|F_n^A|=f_A n+o(n),\quad\mbox{for all $A\in 2^{V(\B f)}$.}$$
 (Of course, the latter condition is automatically true for all $A\in
2^{V(\B f)}$ with $f_A=0$.) Clearly, $e(F_n)=(e(\B f)+o(1))n$, where
$e(\B f)=\sum_{A\in 2^{V(\B f)}} f_A\,|A|$, so we call $e(\B f)$
the {\em size} of $\B f$. Also, the {\em order} of $\B f$ is $v(\B
f)=|V(\B f)|$ and the {\em degree} of $x\in V(\B f)$ is
 $$d(x)=\sum_{A\in 2^{V(\B f)}\atop A\ni x} f_A.$$
 Clearly, $e(\B f)=\sum_{x\in V(\B f)} d(x)$.

For example, given $t\in\I R_{>0}$, the sequence $(K_{s,\lceil t
n\rceil})_{n\in\I N}$ is the dilatation of $\B k_{s,t}$, where the
symbol $\B k_{s,t}$ will be reserved for the weight on $[s]$ which has
value $t$ on $[s]$ and zero otherwise. (We assume that
$V_1(K_{s,\lceil t n\rceil})=[s]$.) It is not hard to see that any
sequence of bipartite graphs described in the abstract is in fact a
dilatation of some weight.

We write $F\subset \B f$ if for some bipartition $V(F)=V_1(F)\cup
V_2(F)$ there is an injection $h:V_1(F)\to V(\B f)$ such that for any
$A\subset V_1(F)$ dominated by a vertex of $V_2(F)$ there is $B\subset
V(\B f)$ with $h(A)\subset B$ and $f_B>0$. This notation is justified
by the following trivial lemma.

\begin{lemma}\label{lm:\pl:Gf}
 Let $(F_n)_{n\in\I N}$ be a dilatation of $\B f$. If $F\subset \B f$,
then $F$ is a subgraph of $F_n$ for all sufficiently large
$n$. Otherwise, which is denoted by $F\not\subset \B f$, no $F_n$
contains $F$.\qed\end{lemma}

Let $\B f$ and $\B g$ be weights. Assume that $v(\B f)\le v(\B g)$ by
adding new vertices to $V(\B g)$ and letting $\B g$ be zero on all new
sets. We write $\B f\subset \B g$ if there is an injection $h:V(\B
f)\to V(\B g)$ and numbers $(w_{AB}\ge 0)_{A\in 2^{V(\B f)},\, B\in
2^{V(\B g)}}$ such that
 \begin{eqnarray*}
 \forall A\in 2^{V(\B f)},\ \forall B\in 2^{V(\B g)}&& h(A)\not\subset
B\imply w_{AB}=0,\\
 \forall A\in 2^{V(\B f)}&&\sum_{B\in 2^{V(\B g)}\atop B\supset h(A)}
w_{AB}\ge f_A,\\
 \forall B\in 2^{V(\B g)}&& \sum_{A\in 2^{V(\B f)}\atop h(A)\subset B}
w_{AB}\le g_B.\end{eqnarray*}

This can be viewed as a fractional analogue of the subgraph relation
$F\subset G$: $h$ embeds $V_1(F)$ into $V_1(G)$ and $w_{AB}$ says how
much of $F^A\subset V_2(F)$ is mapped into $G^B$. The fractional
$\subset$-relation enjoys many properties of the discrete one. For
example, $d(x)\le d(h(x))$ for any $x\in V(\B f)$:
 \begin{equation}\label{eq:\pl:d(h(x))} d(x)=\sum_{A\in 2^{V(\B f)}\atop A\ni x} f_A \le \sum_{A\in 2^{V(\B
f)}\atop A\ni x} \sum_{B\in 2^{V(\B g)}\atop B\supset h(A)} w_{A,B}
\le  \sum_{B\in 2^{V(\B g)}\atop B\ni h(x)} \sum_{A\in 2^{V(\B
f)}\atop h(A)\subset B} w_{A,B} \le \sum_{B\in 2^{V(\B g)}\atop B\ni
h(x)} g_B=d(h(x)).\end{equation}

The following result is not difficult and, in fact, we will implicitly
prove a sharper version later (with concrete estimates of $\e$), so we
omit the proof.

\begin{lemma}\label{lm:\pl:fg} Let $(F_n)_{n\in\I N}$ and $(G_n)_{n\in\I N}$ be dilatations
of $\B f$ and $\B g$ respectively. Then $\B f\subset \B g$ implies
that for any $\e>0$ there is $n_0$ such that $F_n\subset G_m$ for any
$n\ge n_0$ and $m\ge (1+\e)n$. Otherwise, which is denoted by $\B
f\not\subset \B g$, there is $\e>0$ and $n_0$ such that
$F_n\not\subset G_m$ for any $n\ge n_0$ and $m\le (1+\e)n$.\qed\end{lemma}

An {\em $r$-colouring} $\B c$ of
$\B g$ is a sequence $(c_{A_1,\dots,A_r})$ of non-negative reals
indexed by $r$-tuples of disjoint subsets of $V(\B g)$ such that
 \begin{equation}\label{eq:\pl:c} \sum_{A_1\cup\cdots\cup A_r=A} c_{A_1,\dots,A_r} > g_A,\quad
\mbox{for all $A\in 2^{V(\B g)}$.}\end{equation}
 The {\em $i$-th colour subweight} $\B c_i$ is defined by $V(\B
c_i)=V(\B g)$ and
 \begin{equation}\label{eq:\pl:ci} c_{i,A}= \sum_{A_1,\dots,A_r\atop A_i=A}
c_{A_1,\dots,A_r},\quad A\in 2^{V(\B g)}.\end{equation}
 The analogy: to define an $r$-colouring of $G$, it is enough to
define, for all disjoint $A_1,\dots,A_r\subset V_1(G)$, how many
vertices of $G^{A_1\cup\cdots\cup A_r}$ are connected, for all
$i\in[r]$, by colour $i$ precisely to $A_i$. Following this analogy,
there should have been the equality sign in~\eqref{eq:\pl:c}; however, the
chosen definition will make our calculations less messy later.

\section{Existence of Limit}

Let $r\ge q\ge1$. Consider a sequence $\B F=(\B F_1,\dots,\B F_r)$,
where $\B F_i=\B f_i$ is a weight for $i\in[q]$ and $\B F_i=F_i$ is a
bipartite graph for $i\in[q+1,r]$. Assume that $\B F_i$ does not have
an {\em isolated vertex} (that is, $x\in V(\B F_i)$ with $d(x)=0$),
$i\in[r]$. We say that a weight $\B g$ {\em arrows} $\B F$ (denoted by
$\B g\to\B F$) if for any $r$-colouring $\B c$ of $\B g$ we have
$\B F_i\subset \B c_i$ for some $i\in[r]$. Define
 \begin{equation}\label{eq:\pl:Bsr}\hat r(\B F)=\inf\{e(\B g)\mid \B g\to \B F\}.\end{equation}

The definition~\eqref{eq:\pl:Bsr} imitates that of the size Ramsey number and
we will show that these are very closely related indeed. However, we
need a few more preliminaries.

Observe that $\hat r(\B F)<\infty$ by considering $\B k_{a,b}$ which
arrows $\B F$ if, for example, $a=1-r+\sum_{i=1}^r v(\B F_i)$ and $b$
is sufficiently large, cf.~\eqref{eq:\pl:On}. Let $l$ be an integer greater
than $\hat r(\B F)/\dmin$, where $\dmin=\sum_{i=1}^q d_i$ and
 $$d_i=\min\{d_{\B f_i}(x)\mid x\in V(\B f_i)\}>0,\quad i\in[q].$$

\begin{lemma}\label{lm:\pl:l} $\hat r(\B F)=\hat r_l(\B F)$, where $r_l(\B F)=\min\{e(\B
g)\mid \B g\to \B F,\ v(\B g)\le l\}$.\end{lemma}
 \smallskip{\it Proof.}  Let $\e>0$ be any real smaller than $\dmin$. Let $\B g\to\B F$ be
a weight with $v(\B g)>l$ and $e(\B g)\le \hat r(\B F)+\e$. To prove
the theorem, it is enough to construct $\B g'\to\B F$ with $e(\B
g')\le e(\B g)$ and $v(\B g')=v(\B g)-1$.

We have $d(x)\le (\hat r(\B F)+\e)/(l+1)< \dmin$ for some $x\in V(\B
g)$. Choose $\delta>0$ with $\delta+d_i d(x)/\dmin <d_i$ for any
$i\in[q]$.  Define the weight $\B g'$ on $V(\B g)\setminus\{x\}$ by
$g_A'=g_A+g_{A\cup\{x\}}$, $A\in2^{V(\B g')}$. Clearly, $e(\B g')=e(\B
g)-d(x)\le e(\B g)$.

We claim that $\B g'$ arrows $\B F$. Suppose that this is not true and
let $\B c'$ be an $\B F$-free $r$-colouring of $\B g'$. We can assume
that
 $$ \sum_{A_1\cup\cdots\cup A_r=A}c_{A_1,\dots,A_r}' \le g_A'
+\delta/r^l,\quad \mbox{for any $A\in 2^{V(\B g')}$.}$$
 Define $\B c$ by
 $$c_{A_1,\dots,A_r}= \left\{\begin{array}{ll}
 \frac{\lambda_{A\setminus\{x\}} d_i}{\dmin} \cdot
c_{A_1,\dots,A_{i-1},A_i\setminus\{x\},A_{i+1},\dots,A_r}',&
x\in A_i,\ i\in[q],\\
 0,& x\in A_{q+1}\cup\dots\cup A_r,\\
 (1-\lambda_{A}) \cdot c_{A_1,\dots,A_r}',& x\not\in
A,\end{array}\right.$$
 where we denote $A=A_1\cup\cdots\cup A_r$,
$\lambda_A=g_{A\cup\{x\}}/g_A'$ if $g_A'>0$, and $\lambda_A=1/2$ if
$g_A'=0$. The reader can check that $\B c$ is an $r$-colouring of $\B
g$.

By the assumption on $\B g$, we have $\B F_i\subset \B c_i$ for some
$i\in[r]$. But this embedding cannot use $x$ because for $i\in[q+1,r]$
we have $d_{\B c_i}(x)=0$ while for $i\in[q]$
 \begin{eqnarray*}
 d_{\B c_i}(x) &=& \sum_{A_1,\dots,A_r\subset V(\B g')}
c_{A_1,\dots,A_{i-1},A_i\cup\{x\},A_{i+1},\dots,A_r} \ =\ \sum_{A\in
2^{V(\B g')}} \frac{\lambda_A d_i}{\dmin} \sum_{A_1\cup\cdots\cup
A_r=A} c_{A_1,\dots,A_r}'\\ & \le& \sum_{A\in 2^{V(\B g')}}
\frac{\lambda_A d_i}{d_0} (g_A' +\delta/r^l)\ \le\
\frac{d_i\delta}{\dmin} + \frac{d_i}{\dmin} \sum_{A\in 2^{V(\B g')}}
g_{A\cup\{x\}}\ \le\ \delta + d_i
\frac{d(x)}{\dmin}\ <\ d_i\end{eqnarray*}
 is too small, see~\eqref{eq:\pl:d(h(x))}. But $c_{i,A}\le c_{i,A}'$ for $A\in
2^{V(\B g')}$; hence, $\B F_i\subset \B c_i'$, which is the desired
contradiction.\qed \medskip

Hence, to compute $\hat r(\B F)$ it is enough to consider $\B
F$-arrowing weights on $L=[l]$ only.

\begin{lemma}\label{lm:\pl:minimum} There exists $\B g\to \B F$ with $V(\B g)\subset L$ and
$e(\B g)=\hat r(\B F)$. (And we call such a weight {\em
extremal}.)\end{lemma}
 \smallskip{\it Proof.}  Let $\B g_n\to \B F$ be a sequence with $V(\B g_n)\subset L$
such that $e(\B g_n)$ approaches $\hat r(\B F)$. By choosing a
subsequence, assume that $V(\B g_n)$ is constant and
$g_A=\lim_{n\to\infty} g_{n,A}$ exists for each $A\in2^L$. Clearly,
$e(\B g)=\hat r(\B F)$ so it remains to show that $\B g\to\B F$.

Let $\B c$ be an $r$-colouring of $\B g$. Let $\delta$ be the smallest
slack in inequalities~\eqref{eq:\pl:c}. Choose sufficiently large $n$ so
that $|g_{n,A}-g_A|< \delta$ for all $A\in 2^L$. We have 
 $$\sum_{A_1\cup\cdots\cup A_r=A} g_{A_1,\dots,A_r} \ge g_A+\delta
>g_{n,A},\quad A\in 2^L,$$
 that is, $\B c$ is a colouring of $\B g_n$ as well. Hence,
$\B F_i\subset \B c_i$ for some $i$, as required.\qed \medskip

Now we are ready to prove our general theorem. Its proof essentially
takes care of itself. We just exploit the parallels between the
fractional and discrete universes, which, unfortunately, requires
messing around with various constants.

\begin{theorem}\label{th:\pl:general} Let $(F_{i,n})_{n\in\I N}$ be a dilatation of $\B f_i$,
$i\in[q]$. Then, for all sufficiently large $n$,
 \begin{equation}\label{eq:\pl:general} \hat r(\B F)n-M(1+\fmax)\le \hat
r(F_{1,n},\dots,F_{q,n},F_{q+1},\dots,F_r)\le\hat r(\B
F)n+M(1+\fmax),\end{equation}
 where $\fmax=\max\{\,|\,|F_i^A| - f_{i,A}n\,|\mid i\in[q], A\in V(\B
f_i)\}$ and $M=M(\B F)$ is some constant.

In particular, the limit $\lim_{n\to\infty}\hat
r(F_{1,n},\dots,F_{q,n},F_{q+1},\dots,F_r)/n$ exists.\end{theorem}
 \smallskip{\it Proof.}  Let $\vmax=\max\{v(F_i)\mid i\in[r]\}$,
$\cone=2^{\vmax}(\fmax+1)$, and $\ctwo=r^l\cone+1$.

We will prove that
 \begin{equation}\label{eq:\pl:upper} \hat r(F_1,\dots,F_r)\le \hat r(\B
F)n+2^ll(\ctwo+1),\quad n\ge 1.\end{equation}

By Lemma~\ref{lm:\pl:minimum} choose an extremal weight $\B g$ on $L$. Define
a bipartite graph $G$ as follows. Choose disjoint from each other (and
from $L$) sets $G^A$ with $|G^A| = \lceil g_An + \ctwo \rceil$, $A\in
2^L$. Let $V(G)=L\cup (\cup_{A\in 2^L} G^A)$. In $G$ we connect $x\in
L$ to all of $G^A$ if $x\in A$. These are all the edges. Clearly,
 $$e(G)=\sum_{A\in 2^L} |G^A|\, |A| \le 2^ll(\ctwo +1)+\sum_{A\in 2^L}
g_An\,|A| \le \hat r(\B F)n + 2^ll(\ctwo +1),$$
 as required. Hence, it is enough to show that $G$ has the arrowing
property.

Consider any $r$-colouring $c:E(G)\to[r]$. For disjoint sets
$B_1,\dots,B_r\subset L$, let
 \begin{eqnarray*}
 C_{B_1,\dots,B_r}&=&\{y\in G^B\mid \forall i\in[r],\ \forall x\in B_i\
c(\{x,y\})=i\},\\
 c_{B_1,\dots,B_r}&=&\left\{\begin{array}{ll}
 (|C_{B_1,\dots,B_r}|-\cone)/n,&\mbox{if $|C_{B_1,\dots,B_r}|\ge \cone$},\\
 0,& \mbox{otherwise,}\end{array}\right.\end{eqnarray*}
 where $B=B_1\cup\cdots\cup B_r$. Clearly, $nc_{B_1,\dots,B_r}\ge
|C_{B_1,\dots,B_r}|-\cone$; hence,
 $$n\sum_{B_1\cup\cdots\cup B_r=B} c_{B_1,\dots,B_r} \ge -r^{|B|}
\cone+\sum_{B_1\cup\cdots\cup B_r=B} |C_{B_1,\dots,B_r}|\ge -r^{l} \cone +
|G^B|> g_B,$$
 that is, $\B c$ is an $r$-colouring of $\B g$. Hence, $\B F_i\subset
\B c_i$ for some $i\in[r]$.

Suppose that $i\in[q]$. By the definition, we find appropriate $h:V(\B
f_i)\to L$ and $\B w$. We aim at proving that $F_{i,n}\subset G_i$,
where $G_i\subset G$ is the colour-$i$ subgraph. Partition
$F_{i,n}^A=\cup_{B\supset h(A)} W_{A,B}$ so that $W_{A,B}=\emptyset$ if
$w_{A,B}=0$ and $|W_{A,B}|\le \lfloor w_{A,B}n+\fmax+1\rfloor$, $A\in
2^{V(\B f_i)}$, $B\in 2^L$. This is
possible for any $A$: if $w_{A,B}=0$ for all $B\in 2^L$ with
$h(A)\subset B$, then $f_{i,A}=0$ and $F^A=\emptyset$; if $w_{A,B}>0$
for at least one $B$, then
 $$\sum_{B\in 2^L\atop w_{A,B}>0} (w_{A,B}n+\fmax)\ge
\fmax+n\sum_{B\in 2^L\atop w_{A,B}>0} w_{A,B} \ge \fmax+f_{i,A}n\ge
|F_{i,n}^A|.$$

Let $B\in 2^L$. If $c_{i,B}=0$, then $|W_{A,B}|=w_{A,B}=0$ for all
$A\in 2^{V(\B f_i)}$. Otherwise,
 $$c_{i,B}n=n\sum_{B_1,\dots,B_r\atop B_i=B}
c_{B_1,\dots,B_i} \le -\cone+\sum_{B_1,\dots,B_r\atop
B_i=B}|C_{B_1,\dots,B_i}|=|G_i^B|-\cone,$$
 and we have
 $$\sum_{A\in 2^{V(F_{i,n})}\atop h(A)\subset B}
|W_{A,B}| \le \sum_{A\in 2^{V(F_{i,n})}\atop h(A)\subset B}
(w_{A,B}n+\fmax+1)\le c_{i,B} n+2^{\vmax}(\fmax+1)\le |G_i^B|.$$
 Hence, we can extend $h$ to the whole of $V(F_{i,n})$ by mapping
$\cup_{h(A)\subset B}W_{A,B}$ injectively into~$B$.

Suppose that $i\in[q+1,r]$. The relation $F_i\subset \B c_i$ means
that there exist appropriate $V_1(F_i)\cup V_2(F_i)=V(F_i)$ and
$h:V_1(F_i)\to L$. We view $h$ as a partial embedding of $F_i$ into
$G_i$ and will extend $h$ to the whole of $V(F_i)$.

Take consecutively $y\in V_2(F_i)$. There is $B_i\subset L$ such that
$c_{i,B_i}>0$ and $h(\Gamma(y))\subset B_i$. The inequality
$c_{i,B_i}>0$ implies that there are disjoint $B_j$'s,
$j\in[r]\setminus\{i\}$, such that $c_{B_1,\dots,B_r}>0$. Each vertex
in $C_{B_1,\dots,B_r}$ is connected by colour $i$ to the whole of
$B_i\supset h(\Gamma(y))$. The inequality $c_{B_1,\dots,B_r}>0$ means
that $|C_{B_1,\dots,B_r}|\ge \cone\ge v(F_i)$, so we can always extend
$h$ to $y$. Hence, we find an $F_i$-subgraph of colour $i$ in this
case.

Thus the constructed graph $G$ has the desired arrowing property, which
proves the upper bound.

Let $\cthree=\max(r^{4l/\dminmin}, 2^{\vmax}\fmax)$, where
$\dminmin=\min_{i\in[q]} \min _{x\in V(\B f_i)} d_{f_i}(x)>0$. As the
lower bound, we show that, for all sufficiently large $n$,
 \begin{equation}\label{eq:\pl:lower} \hat r(F_1,\dots,F_r)\ge \hat r(\B
F)n-4l\cdot2^{4l/\dminmin}\cdot \cthree/\dminmin.\end{equation}

Suppose on the contrary that we can find an arrowing graph $G$
contradicting~\eqref{eq:\pl:lower}. Let $L\subset V(G)$ be the set of vertices
of degree at least $\dminmin n/2$ in $G$. From $\dminmin n|L|/4< e(G)<
l n$ it follows that $|L|\le 4l/\dminmin$. For $A\in 2^L$, define
$g_{A}=(|G^A|+\cthree)/n$.

We have 
 $$\sum_{A\in 2^L} g_A\, |A|\le \frac{2^{4l/\dminmin}\cthree}{n} \cdot
\frac{4l}{\dminmin} +\frac1{n}\sum_{A\in 2^L} |G^A|\, |A|\le
\frac{4l\cdot 2^{4l/\dminmin}\cthree/\dminmin+e(G)}{n} < \hat r(\B
F).$$

Thus there is an $\B F$-free $r$-colouring $\B c$ of $\B g$. We are
going to exhibit a contradictory $r$-colouring of $E(G)$.

For each $B\in 2^L$ choose any disjoint sets $C_{B_1,\dots,B_r}\subset
G^B$ (indexed by $r$-tuples of disjoint sets partitioning $B$) such
that they partition $G^B$ and $|C_{B_1,\dots,B_r}|\le\lfloor
c_{B_1,\dots,B_r}\cdot n\rfloor$. This is possible because
 $$ \sum_{B_1\cup\cdots\cup B_r=B}\lfloor
c_{B_1,\dots,B_r}\cdot n\rfloor \ge g_Bn-r^{4l/\dminmin} \ge |G^B|.$$
 For $j\in[r]$, $x\in B_j$ and $y\in C_{B_1,\dots,B_r}$, colour the
edge $\{x,y\}$ by colour $j$. All the remaining edges of $G$ (namely,
those lying inside $L$ or inside $V(G)\setminus L$) are coloured with
colour $1$.

There is $i\in[r]$ such that $G_i\subset G$, the colour-$i$ subgraph,
contains a forbidden subgraph.

Suppose that $i\in[q]$. Let $h:F_{i,n}\to G_i$ be an embedding. If $n$
is large, then $d(x)>4l/\dminmin+\dminmin n/2$ for all $x\in
V_1(F_{i,n})$, which implies that $h(V_1(F_{i,n}))\subset L$.  Define
for $A\in 2^{V(\B f_i)}$ and $B\in 2^L$ with $B\supset h(A)$ and
$f_{i,A}\not=0$
 $$w_{A,B}=\frac{|h^{-1}(G^B)\cap F_i^A|+\fmax}n.$$
 All other $w_{A,B}$'s are set to zero. For $A\in 2^{V(\B f_i)}$ with
$f_{i,A}\not=0$, we have
 $$\sum_{B\in2^L\atop B\supset h(A)} w_{A,B} \ge
(|F_i^A|+\fmax)/n\ge f_{i,A}.$$
 For $B\in 2^L$ we have
 $$\sum_{A\in 2^{V(\B f_i)}\atop h(A)\subset B} w_{A,B} \le
\frac{2^{\vmax}\fmax}n+\sum_{A\in 2^{V(\B f_i)}\atop h(A)\subset B}
\frac{|h^{-1}(G^B)\cap F_i^A|}n\le
\frac{2^{\vmax}\fmax}n +\frac{|G^B|}n \le g_B,$$
 that is, $h$ (when restricted to $V(\B f_i)$) and $\B w$ demonstrate
that $\B f_i\subset\B c_i$, which is a contradiction.

Suppose that $i\in[q+1,r]$. Let $V_1(F_i)$ consists of those vertices
which are mapped by $h:F_i\to G_i$ into $L$ and let
$V_2(F_i)=V(F_i)\setminus V_1(F_i)$. This is a legitimate bipartition
of $F_i$ because any colour-$i$ edge of $G$ connects $L$ to
$V(G)\setminus L$. Let $y\in V_2(F_i)$. The sets $C_{B_1,\dots,B_r}$'s
partition $V(G)\setminus L$; let $y\in C_{B_1,\dots,B_r}$. If
$\{y,z\}\in E(F_i)$, $z\in V_1(F_i)$, then $B_i\ni h(z)$, which,
together with $c_{B_1,\dots,B_r}>0$ shows that $F_i\subset \B
g_i$. This contradiction proves the theorem.\qed \medskip

\section{Complete Bipartite Graphs}

Here we will compute asymptotically the size Ramsey number if each
forbidden graph is a complete bipartite graph. More precisely, we
show that in order to do this it is enough to consider only complete
bipartite graphs having the arrowing property.

\begin{theorem}\label{th:\pl:mcomplete} Let $r\ge2$ and $q\ge1$. Suppose that we are given
$t_1,\dots,t_q\in \I R_{>0}$ and $s_1,\dots,s_r,t_{q+1},\dots,t_r\in\I
N$ such that $t_i\ge s_i$ for $i\in[q+1,r]$. Then there exist $s\in\I
N$ and $t\in\I R_{>0}$ such that $\B k_{s,t}\to \B F$ and $\hat r(\B
F)=e(\B k_{s,t})=st$, where
 $$\B F=(\B k_{s_1,t_1}, \dots, \B k_{s_q,t_q}, K_{s_{q+1},t_{q+1}},
\dots, K_{s_r,t_r}).$$\end{theorem}
 \smallskip{\it Proof.}  Let us first describe an algorithm finding extremal $s$ and
$t$. Some by-product information gathered by our algorithm will be
used in the proof of the extremality of $\B k_{s,t}\to \B F$.

Choose $l\in\I N$ bigger than $\hat r(\B F)/\tmin$, where
$\tmin=\sum_{i=1}^q t_i$, which is the same definition of $l$ as that
before Lemma~\ref{lm:\pl:l}.

We claim that $l>\sigma$, where $\sigma=\sum_{i=1}^r(s_i-1)$. Indeed,
take any extremal $\B f\to\B F$ without isolated vertices. The proof
of Lemma~\ref{lm:\pl:l} implies that $d(x)\ge \tmin$ for any $x\in V(\B
f)$. Note that necessarily $v(\B f)>\sigma$, which implies the claim.

For each integer $s\in[\sigma+1,l]$, let $t_s'>0$ be the infimum of
$t\in\I R$ such that $\B k_{s,t}\to \B F$.  Also, let $\R_s$ be the
set of all sequences $\B a=(a_1,\dots,a_r)$ of non-negative integers
with $a_i=s_i-1$ for $i\in[q+1,r]$ and $\sum_{i=1}^r a_i=s$. For a
sequence $\B a=(a_1,\dots,a_r)$ and a set $A$ of size $\sum_{i=1}^r
a_i$, let $\Rbinom{A}{\B a}$ consist of all sequences $\B
A=(A_1,\dots,A_r)$ of sets partitioning $A$ with $|A_i|=a_i$,
$i\in[r]$.

We claim that $t_s'$ is $\sol(L_s)$, the extremal value of the
following linear program $L_s$: ``\sl Find $\sol(L_s)=\max\sum_{\B
a\in\R_s} w_{\B a}$ over all sequences $(w_{\B a})_{\B a\in\R_s}$ of
non-negative reals such that\begin{equation}\label{eq:\pl:mb1}
 \sum_{\B a\in\R_s} w_{\B a} \binom {a_i}{s_i}\le t_i\binom s{s_i},
\quad \mbox{for all $i\in[q]$.''}\end{equation}\rm

\claim1 The weight $\B k_{s,t}$ does not arrow $\B F$
for $t<\sol(L_s)$. 

\smallskip\noindent To prove this, let
 $$\lambda=\frac{t+\sol(L_s)}{2\sol(L_s)}<1\quad \mbox{and}\quad
\e=\frac{1-\lambda}{2 r^l}\, \min\{t_i\mid i\in[q]\}>0.$$
 Let $V(\B k_{s,t})=[s]$. Define an $r$-colouring $\B c$ of $\B
k_{s,t}$ by
 $$c_{\B A}=\frac{\lambda w_{|A_1|,\dots,|A_r|}}{\binom
s{|A_1|,\dots,|A_r|}},\quad \B a\in\R_s,\ \B A\in\Rbinom{[s]}{\B a},$$
 $c_{B,\emptyset,\dots,\emptyset} =\e$, $B\subsetneq [s]$,  
while all other $c$'s are zero. It is indeed a colouring: 
 $$\sum_{\B a\in\R_s} \sum_{\B A\in\Rbinom{[s]}{\B a}} c_{\B A} =
\sum_{\B a\in\R_s} \lambda w_{\B a} =\lambda\, \sol(L_s) >t.$$
 We have $\B k_{s_i,t_i}\not\subset\B c_i$ for $i\in[q]$. For example,
for $i=1$ and any $ S\in\binom{[s]}{s_1}$, we have 
 $$\sum_{\B a\in\R_s} \sum_{\B A\in\Rbinom{[s]}{\B a}\atop A_1\supset
S} c_{\B A} =\sum_{\B a\in\R_s\atop a_1\ge s_1}
\frac{\binom{s-s_1}{a_1-s_1,a_2,\dots,a_r} \lambda w_{\B a}}{\binom
s{a_1,\dots,a_r}} =\lambda\sum_{\B a\in\R_s}\frac{\binom {a_1}{s_1}
w_{\B a}}{\binom s{s_1}}\le \lambda t_1 <t_1-\sum_{B\subsetneq
[s]\atop B\supset S} c_{B,\emptyset,\dots,\emptyset}.$$

Also, $K_{s_i,t_i}\not\subset \B c_i$ for $i\in[q+1,r]$ because
$c_{A_1,\dots,A_r}=0$ whenever $|A_i|\ge s_i$ for some
$i\in[q+1,r]$. Claim~1 is proved.

\claim2 $\B k_{s,t}\to \B F$ for any $t> \sol(L_s)$.

\smallskip\noindent Suppose that the claim is not true and we can find
an $\B F$-free $r$-colouring $\B c$ of $\B k_{s,t}$. By the
definition, $c_{A_1,\dots,A_r}=0$ whenever $|A_i|\ge s_i$ for some
$i\in[q+1,r]$. If some $c_{A_1,\dots,A_r}=c>0$ with $|A_i|\le s_i-2$
for some $i\in[q+1,r]$, then $A_j\not=\emptyset$ for some $j\in[q]$,
so we can pick $x\in A_j$ and set $c_{A_1,\dots,A_r}=0$ while
increasing $c_{\dots,A_j\setminus\{x\},\dots,A_i\cup\{x\},\dots}$ by
$c$.  Clearly, $\B c$ remains an $\B F$-free colouring. Thus, we can
assume that all the $c$'s are zero except those of the form $c_{\B
A}$, $\B A\in\Rbinom{[s]}{\B a}$ for some $\B a\in\R_s$. Now,
retracing back our proof of Claim~1, we obtain a feasible solution
$w_{\B a}=\sum_{\B A\in\Rbinom{[s]}{\B a}} c_{\B A}$, $\B a\in\R_s$, to
$L_s$ with a larger objective function, which is a contradiction. The
claim is proved.\smallskip

Thus, $t_s'=\sol(L_s)$ and $m_u=\min\{st_s'\mid s\in[\sigma+1,l]\}$ is
an upper bound on $\hat r(\B F)$. Let us show that in fact $\hat r(\B
F)=m_u$.

We rewrite the definition of $\hat r(\B F)$ so that we can apply the
Farkas Lemma. The proof of the following easy claim is left to the
reader.

\claim3 $\hat r(\B F)=\inf e(\B g)$ over all weights $\B g$ on $L$
such that there do not exist non-negative reals $(c_{\B A})_{\B A\in
\Rbinom{A}{\B a},\, \B a\in\R_{|A|},\, A\in 2^L}$ with the following
properties\begin{eqnarray*}
 \sum_{\B a\in\R_{|A|}} \sum_{\B A\in\Rbinom{A}{\B a}} c_{\B A}& \ge&
g_A,\quad A\in 2^L, \\
 \sum_{A\in 2^L} \sum_{\B a\in\R_{|A|}} \sum_{\B A\in\Rbinom{A}{\B a}
\atop A_i\supset S} c_{\B A} &\le& t_i,\quad i\in[q],\ S\in\binom
L{s_i}.  \end{eqnarray*}

Let $\B g$ be any feasible solution to the above problem. By the Farkas
Lemma there exist $x_A\ge 0$, $A\in2^L$, and $y_{i,S}\ge 0$,
$i\in[q]$, $S\in\binom L{s_i}$, such that\begin{eqnarray}
  \sum_{i=1}^q \sum_{S\in\binom {A_i}{s_i}} y_{i,S}&\ge& x_A,\quad
A\in2^L,\ \B a\in\R_{|A|},\ \B A\in\Rbinom{A}{\B a}, \label{eq:\pl:mxyz1}\\
 \sum_{i=1}^q t_i \sum_{S\in \binom L{s_i}} y_{i,S} &<& \sum_{A\in2^L}
g_Ax_A.\label{eq:\pl:mxyz2}\end{eqnarray}

We deduce that $x_A\le 0$ (and hence $x_A=0$) if $|A|\le\sigma$ by
considering~\eqref{eq:\pl:mxyz1} for some $\B A$ with $|A_i|\le s_i-1$, each
$i\in[r]$.

For each $A$ with $a:=|A|>\sigma$ repeat the following. Let $(w_{\B
a})_{\B a\in\R_a}$ be an extremal solution to $L_{a}$.  For each $\B
a\in\R_{a}$, take the average of~\eqref{eq:\pl:mxyz1} over all $\B
A\in\Rbinom{A}{\B a}$, multiply it by $w_{\B a}$, and add all these
equalities together to obtain the following.\begin{eqnarray*}
 x_A t_a' &=& \sum_{\B a\in\R_{a}}\sum_{\B A\in\Rbinom{A}{\B a}}
w_{\B A} x_A\ \le\ \sum_{\B a\in\R_{a}} \frac{w_{\B a}}{\binom
a{a_1,\dots,a_r}} \sum_{\B A\in\Rbinom{A}{\B a}}
 \sum_{i=1}^q\sum_{S\in \binom{A_i}{s_i}} y_{i,S}\\
 &=& \sum_{i=1}^q \sum_{S\in \binom {A}{s_i}} y_{i,S}
\sum_{\B a\in\R_a\atop a_i\ge s_i} \frac{w_{\B
a}\binom{a-s_i}{a_1,\dots,a_{i-1},a_i-s_i,a_{i+1},\dots,a_q}}{\binom
a{a_1,\dots,a_r}}\\
 &=& \sum_{i=1}^q \sum_{S\in \binom A{s_i}} y_{i,S}
\sum_{\B a\in\R_{a}\atop a_i\ge s_i} \frac{w_{\B
a}\binom{a_i}{s_i}}{\binom a{s_i}}\ \le\ \sum_{i=1}^q t_i\sum_{S\in
\binom A{s_i}} y_{i,S}.
 \end{eqnarray*}
 (In the last inequality we used~\eqref{eq:\pl:mb1}.) 

Substituting the obtained inequalities on the $x_A$'s into~\eqref{eq:\pl:mxyz2}
we obtain
 $$ \sum_{i=1}^q t_i \sum_{S\in \binom L{s_i}} y_{i,S}<
\sum_{A\in2^L\atop|A|>\sigma} \frac{g_A}{t_s'} \sum_{i=1}^q
t_i\sum_{S\in\binom A{s_i}} y_{i,S}.$$
 As the $y_{i,S}$'s are non-negative, some of these variables has a
larger coefficient on the right-hand side. Let it be
$y_{i,S}$. We have
 \begin{equation}\label{eq:\pl:ti}t_i<t_i\sum_{A\in\binom{L}{>\sigma}\atop A\supset S}
\frac{g_A}{t_{|A|}'}\le \frac{t_i}{m_u}\sum_{A\in 2^L} g_A|A|.\end{equation}
 The last inequality follows from the fact that for any integer $a>
\sigma$, we have $1/t_a'\le a/m_u$, which in turn follows from the
definition of $m_u$. Hence, $e(\B g)\ge m_u$ as required.\qed \medskip

\begin{corollary}\label{cr:\pl:complete} Let $r\ge q\ge1$, $t_1,\dots,t_q\in \I R_{>0}$ and
$s_1,\dots,s_r,t_{q+1},\dots,t_r\in\I N$ such that $t_i\ge s_i$ for
$i\in[q+1,r]$. For $i\in[q]$, let $(t_{i,n})_{n\in \I N}$ be an
integer sequence with $t_{i,n}=t_in+o(n)$. Define
 $$\B F_n=(K_{s_1,t_{1,n}},\dots, K_{s_q,t_{q,n}},
K_{s_{q+1},t_{q+1}}, \dots, K_{s_r,t_r}).$$
 Let $l\in\I N$ be larger than $\lim_{n\to\infty} \hat r(\B F_n)/\tmin
n$, where $\tmin=\sum_{i=1}^q t_i$. Then
 \begin{equation}\label{eq:\pl:complete}\lim_{n\to\infty}\frac{\hat r(\B F_n)}n =\lim_{n\to\infty}
\frac{\min\{e(K_{s,t})\mid s\le l,\ K_{s,t}\to\B F_n\}}n.\qed\end{equation}\end{corollary}

In other words, in order to compute the limit in
Corollary~\rcr{complete} it is sufficient to consider only complete
bipartite graphs arrowing $\B F_n$. It seems that there is no simple
general formula, but the proof of Theorem~\ref{th:\pl:mcomplete} gives
an algorithm for computing $\hat r(\B F)$. The author has realized the
algorithm as a C program calling the {\tt lp\_solve} library. (The
latter is a freely available linear programming software, currently
maintained by Michel Berkelaar~\cite{lpsolve}). Later, David Avis
rewrote the program to be linked with his {\tt lrslib}
library~\cite{lrs}. The latter library has the advantage that its
arithmetic is exact (whilst {\tt lp\_solve} operates with reals), so
that any computed limit can be considered as proved (provided the
realisation of our algorithm is correct). The reader is welcome to
experiment with the program; its source is included at the end of this
paper.

For certain series of parameters we can get a more explicit
expression. First, let us treat the case $q=1$, that is, when only the
first forbidden graph dilates with $n$.  We can assume that $t_1=1$ by
scaling $n$.

\begin{theorem}\label{th:\pl:q=1} Let $q=1$ and $r\ge 2$. Then for any
$s_1,\dots,s_r,t_2,\dots,t_r\in\I N$ with $t_i\ge s_i$, $i\in[2,r]$,
we have
 $$\hat r(K_{s_1,n},K_{s_2,t_2},\dots,K_{s_r,t_r})=n\cdot \min\left\{s\,
\frac{(s)_{s_1}}{(s-s')_{s_1}}\mid s\in\I N_{>\sigma}\right\}
+ O(1),$$
 where $s'=\sigma-s_1+1$, $\sigma=\sum_{i=1}^r (s_i-1)$ and
$(s)_k=s(s-1)\dots(s-k+1)$.\end{theorem}
 \smallskip{\it Proof.}  The Problem~$L_s$ has only one variable
$w_{s-s',s_2-1,\dots,s_r-1}$. Trivially, $t_s'=\binom s
{s_1}/\binom{s-s'}{s_1}= (s)_{s_1}/(s-s')_{s_1}$, and the
theorem follows.\qed \medskip

In the case $s_1=1$ we obtain the following formula (with a little bit
of algebra).

\begin{corollary}\label{cr:\pl:q=1} For any $s_2,\dots,s_r,t_2,\dots,t_r\in\I N$ with $t_i\ge
s_i$, $i\in[2,r]$, we have
 $$\hat r(K_{1,n},K_{s_2,t_2},\dots,K_{s_r,t_r})=4\left
(1-r+\sum_{i=2}^r s_i\right) n +O(1).\qed$$\end{corollary}

Another case with a simple formula for $\hat r(\B F)$ is $q=2$,
$s_1=s_2$, and $t_1=t_2$. Again, without loss of generality we can
assume that $t_1=t_2=1$.

\begin{theorem}\label{th:\pl:q=2} Let $q=2$ and $r\ge 2$. Then for any
$s,s_3,\dots,s_r,t_3,\dots,t_r\in\I N$ with $t_i\ge s_i$, $i\in[3,r]$,
we have
 \begin{equation}\label{eq:\pl:q=2}\hat r(K_{s,n},K_{s,n},K_{s_3,t_3},\dots,K_{s_r,t_r})=n\cdot
\min\left\{a\cdot f(a)\mid a\in\I N_{>\sigma}\right\} + O(1),\end{equation}
 where $\sigma=2s-r+\sum_{i=3}^r s_i$
and
 $$f(a)=\frac{2\binom as}{\binom{\lfloor a'/2\rfloor}s +
\binom{\lceil a'/2\rceil}s},$$
 with $a'=a-\sum_{i=3}^r (s_i-1)$.\end{theorem}
 \smallskip{\it Proof.}  Let $a\in\I N_{>\sigma}$ and let $(w_{\B a})_{\B a\in \R_a}$
be an extremal solution to $L_a$. (Where we obviously define
$s_1=s_2=s$ and $t_1=t_2=1$.) Excluding the constant indices in $w_{\B a}$,
we assume that the index set $\R_a$ consists of pairs of integers
$(a_1,a_2)$ with $a_1+a_2=a'$.

Clearly, $w_{a_1,a_2}'=\frac12(w_{a_1,a_2}+w_{a_2,a_1})$,
$(a_1,a_2)\in\R_a$, is also an extremal solution, so we can assume
that $w_{a_1,a_2}=w_{a_2,a_1}$ for all $(a_1,a_2)\in\R_a$.

If $w_{a_1,a_2}=c>0$ for some $a_1<\lfloor a'/2\rfloor$, then we can
set $w_{a_1,a_2}=w_{a_2,a_1}=0$, while increasing $w_{\lfloor
a'/2\rfloor,\lceil a'/2\rceil}$ and $w_{\lceil a'/2\rceil,\lfloor
a'/2\rfloor}$ by $c$. The easy inequality
 $$\binom{b+1}s+\binom{a'-b-1}s -\binom{b}s - \binom{a'-b}s =\binom
b{s-1} -\binom{a'-b-1}{s-1}<0,\quad s-1\le b< \lfloor a'/2\rfloor$$
 implies inductively that the left-hand side of~\eqref{eq:\pl:mb1} strictly
decreases while the objective function $\sum_{\B a\in\R_a} w_{\B a}$
does not change, which clearly contradicts the minimality of $\B w$.

Now we deduce that, for {\em any\/} extremal solution $(w_{\B a})_{\B
a\in\R_a}$, we have $w_{a_1,a_2}=0$ unless $\{a_1,a_2\}=\{\lfloor
a'/2\rfloor, \lceil a'/2\rceil\}$; moreover, it follows now that
necessarily $w_{\lfloor a'/2\rfloor, \lceil a'/2\rceil}=w_{\lceil
a'/2\rceil,\lfloor a'/2\rfloor}$. Hence, $t_a'=f(a)$, which proves the
theorem.\qed \medskip

The special case $r=2$ of Theorem~\ref{th:\pl:q=2} answers the question of
Erd\H os, Faudree, Rousseau and
Schelp~\cite[Problem~B]{erdos+faudree+rousseau+schelp:78}, who asked
for the asymptotics of $\hat r(K_{s,n},K_{s,n})$. Unfortunately, we do
not think that the formula~\eqref{eq:\pl:q=2} can be simplified further in this
case.

Finally, let us consider the case $(r,s)=(2,2)$ of Theorem~\ref{th:\pl:q=2}
in more detail. It is routine to check that Theorem~\ref{th:\pl:q=2} implies
that $\hat r(K_{2,n},K_{2,n})=18n+O(1)$. But we are able to show
that~\eqref{eq:\pl:frs83} holds for all sufficiently large $n$, which is done
by showing that $(K_{2,n},K_{2,n})$-arrowing graph with $18n+o(n)$
edges can have at most $3$ vertices of degree at least $n$ for all
large $n$.

\begin{theorem}\label{th:\pl:q=2,s=2} There is $n_0$ such that, for all $n>n_0$, we have $\hat
r(K_{2,n},K_{2,n})=18n-15$ and $K_{3,6n-5}$ is the only extremal graph
(up to isolated vertices).\end{theorem}
 \smallskip{\it Proof.}  For $n\in\I N$ let $G_n$ be a minimum
$(K_{2,n},K_{2,n})$-arrowing graph. We know that $e(G_n)\le 18n-15$ so
$l_n=|L_n|\le 18$ for all large $n$, where $L_n=\{x\in V(G_n) \mid
d(x)\ge n\}$; let us assume $L_n\subset[18]$. 

\claim1  $l_n\le 3$ for all sufficiently large $n$.

\smallskip\noindent Suppose on the contrary that we can find an
increasing subsequence $(n_i)_{i\in\I N}$ with $l_{n_i}\ge 4$ for all
$i$. Choosing a further subsequence, assume that $L_{n_i}=L$ does not
depend on $i$ and that $g_A=\lim_{i\to\infty} |G_{n_i}^A|/n_i$ exists
for any $A\in2^L$. The argument of Lemma~\ref{lm:\pl:minimum} shows that the
weight $\B g$ on $L$ arrows $(\B k_{2,1},\B k_{2,1})$. 

We have $e(\B g)=18$. It is routine to check that $at_a'>18$ for any
$a\in[4,18]$. The inequality~\eqref{eq:\pl:ti} implies that, for some
$S=\{x,y\}\subset L$, we have $g_A=0$ whenever $|A|>4$ or
$A\not\supset S$. Let $J$ be the set of those $j\in L$ with
$g_{\{x,y,j\}}>0$. We have $\sum_{j\in J} g_{\{x,y,j\}}=6$.

Consider the $2$-colouring $\B c$ of $\B g$ obtained by
letting $c_{A_1,A_2}=2^{-18}/10$ for all disjoint $A_1,A_1\in 2^L$ except
 $$\begin{array}{lclcll}
 c_{\{x,j\},\{y\}}& =& c_{\{y,j\},\{x\}}& = &c_{\{x\},\{y,i\}}\ =\ 
c_{\{y\},\{x,i\}}\ =\ 0.9,\smallskip\\
 c_{\{x,y\},\{j\}} &=& c_{\{j\},\{x,y\}} &=&
(g_{\{x,y,j\}}-3.5)_+/2, \end{array}\quad\quad j\in J,$$  
 where $f_+=f$ if $f>0$ and $f_+=0$ if $f\le 0$. It is easy to check
that neither $\B c_1$ nor $\B c_2$ contains $\B k_{2,1}$: for example,
$\sum_{A\in 2^L: A\supset\{x,y\}} c_{i,A} < (5-3.5)/2 +
0.1<1$. (Recall that $d_{\B g}(x)\ge 1$ for all $x\in L$.) This
contradiction proves Claim~1.\smallskip

Thus, $|L_n|\le 3$ for all large $n$. By the minimality of $G_n$,
$V(G_n)\setminus L_n$ spans no edge and each $x\in V(G_n)\setminus
L_n$ sends at least $3$ edges to $L_n$. (In particular, $|L_n|=3$.)
Thus, disregarding isolated vertices, $G_n=K_{3,m}$. The relation
$G_n\to (K_{2,n},K_{2,n})$ implies that $m\ge 6n-5$, which proves the
theorem.\qed \medskip

\smallskip\noindent{\bf Remark.}  We do not write an explicit expression for $n_0$, although it
should be possible to extract this from the proof (with more algebraic
work) by using the estimates of Theorem~\ref{th:\pl:general}.\smallskip

\section{Generalizations}\label{generalized}

If all forbidden graphs are the same, then one can generalize the
arrowing property in the following way: a graph $G$ {\em
$(r,s)$-arrows} $F$ if for any $r$-colouring of $E(G)$ there is an
$F$-subgraph that receives less than $s$ colours. Clearly, in the case
$s=2$ we obtain the usual $r$-colour arrowing property
$G\to(F,\dots,F)$.

This property was first studied by Ekeles, Erd\H os and F\"uredi (as
reported in~\cite[Section 9]{erdos:81:cn}); the reader can
consult~\cite{axenovich+furedi+mubayi:00} for references to more
recent results.

Axenovich, F\"uredi and Mubayi~\cite{axenovich+furedi+mubayi:00}
studied the generalized arrowing property for bipartite graphs in the
situation when $F$ and $s$ are fixed, $G=K_{n,n}$, and $r$ grows with
$n$.

We can define $\hat r(F,r,s)$ to be the minimal size of a graph which
$(r,s)$-arrows $F$. Our technique extends to the case when $r$ and $s$
are fixed whilst $F$ grows with $n$ (i.e.,\ is a dilatation). Namely,
it should be possible to show the following.

\smallskip
\sl Let $(F_n)_{n\in\I N}$ be a dilatation of a weight $\B f$ and let
$r\ge s$ be fixed. Then the limit $\lim_{n\to\infty} \hat
r(F_n,r,s)/n$ exists; let us denote it by $\hat r(\B f,r,s)$.

We have $\hat r(\B f,r,s)<\infty$ and, in fact, $\hat r(\B f,r,s)=\min
e(\B g)$ over all weights $\B g$ such that for any $r$-colouring $\B
c$ of $\B g$ there is $S=\{i_1,\dots,i_s\}\in \binom{[r]}s$ such that
$\B c_S\supset \B f$, where
 $$c_{S,A}=\sum_{A_{1},\dots,A_{r}}\ c_{A_1,\dots,A_r},\quad A\in
2^{V(\B g)},$$
 where the sum is taken over all disjoint $A_1,\dots,A_r\in 2^L$ with
$A_{i_1}\cup \dots\cup A_{i_s}=A$.\smallskip\rm

We omit the proof as the complete argument would not be very short and
it is fairly obvious how to proceed.

Also, one can consider the following settings. Let $\C F_i$ be a
family of graphs, $i\in[r]$. We write $G\to(\C F_1,\dots,\C F_r)$ if
for any $r$-colouring of $E(G)$, there is $i\in[r]$ and $F\in\C F_i$
such that we have an $F$-subgraph of colour $i$. The task is to
compute the minimum size of a such $G$. Again, we believe that our
method extends to this case as well. But we do not provide any proof,
so we do not present this as a theorem.

\section*{Acknowledgements}

I am grateful to Martin Henk, Deryk Osthus, and G\"unter Ziegler for
helpful discussions and to David Avis for rewritting and improving my
code (and for his nice {\tt lrslib} library!).

\providecommand{\bysame}{\leavevmode\hbox to3em{\hrulefill}\thinspace}

\newpage
\pagestyle{empty}

\section*{C source code}

\small
\renewcommand{\baselinestretch}{0.84}
\hoffset= -1.5 true cm

\begin{verbatim}

/* lrslib hack by D. Avis */
/* asymptotic computation of size Ramsey function for complete
   bipartite graphs; (C) 2000-1 Oleg Pikhurko;
   distributed under GNU General Public Licence: see http://www.gnu.org 

   This program computes 
   $m=\lim_{n\to\infty} \hat r(K_{s_1,t_{1,n}},\dots, K_{s_q,t_{q,n}},
   K_{s_{q+1},t_{q+1}}, \dots,K_{s_r,t_r})/n$,
   where $t_{i,n}=t_in+o(n)$ with each $t_i$ being a fixed integer.

   Input: 
   q r
   s[1] t[1]
   ...
   s[q] t[q]
   min(s[q+1],t[q+1])
   ...
   min(s[r],t[r])

   Output: the value of m.

   Compiling: the code has to be linked with lrs 4.1 which is
   freely available at http://cgm.cs.mcgill.ca/~avis/C/lrs.html

   Description of the algorithm (and the proof of its correctness) can
   be found in the e-print "Asymptotic Size Ramsey Results for
   Bipartite Graphs" by Oleg Pikhurko at http://www.arXiv.org */


#include <stdio.h>
#include "lrslib.h"

/* maximum number of rows and columns in matrix */

#define MAXR 1000

/* computation of the binomial coefficient */
long binom(long, long);
void Binom(long, long, lrs_mp);

int 
main()
{
  lrs_dic *P;                   /* structure for holding current dictionary and indices */
  lrs_dat *Q;                   /* structure for holding static problem data            */

  long m;                       /* number of constraints in the problem                 */
  long d;                       /* number of variables in the problem                   */
  long *num;                    /* numerators for one constraint                        */
  long *den;                    /* denominatorss for one constraint                     */

  lrs_mp_vector Den;            /* denominators for one constraint                      */

  lrs_mp_matrix Mat;             /* to hold b+Ax>=0                                      */
                                /* row zero holds 0 c_1 ... c_d                         */

  lrs_mp_vector output;         /* holds one line of output; ray,vertex,facet,linearity */

  lrs_mp best_bound_num,best_bound_den;
  lrs_mp mpone, lt, temp;
  

  long q, r, sigma=0;
  long t[MAXR];
  long s[MAXR], a[MAXR];
  long i, j, l, first_lp=1;
  long current_column, ncolumns;

  lrs_alloc_mp(best_bound_num);
  lrs_alloc_mp(best_bound_den);
  lrs_alloc_mp(lt);
  lrs_alloc_mp(mpone);
  lrs_alloc_mp(temp);

  itomp (ONE, mpone);

  printf("\nAsymptotic computation of size Ramsey function for complete");
  printf(" bipartite graphs.\nSee the source code for details.  ");

/* lrs initialization */
  if ( !lrs_init ("\n*lramsey:"))
    return 1;

  lrs_set_digits (200L);   /* fix max number of decimal digits */


/* output holds one line of output from dictionary     */
  output = lrs_alloc_mp_vector (MAXR); 

  printf("\nPlease enter q and r: ");

  scanf("%ld %ld", &q,&r);

  s[0]=0;

  for (i=1; i<=q; i++) {
    printf("Enter s[%ld] and t[%ld]: ", i,i);
    scanf("%ld %ld",&(s[i]),&(t[i]));
    sigma+=s[i]-1;
    t[0]+=t[i];
  }

  for (i=q+1; i<=r; i++) {
    printf("Enter min(s[%ld],t[%ld]): ", i,i);
    scanf("%ld",&s[i]);
    sigma+=s[i]-1;
    s[0]+=s[i]-1; 
  }

  l=sigma+1;

  do {

    /* construct all partitions a[1]+...+a[q]=l-s[0] in lex order */
    /* ncolumns counts the number of partitions */

    ncolumns=binom(l-s[0]+q-1,q-1); 

/* q contraints and ncolumns variables */

    /* temporary arrays for setting up matrix */

    den = calloc ((ncolumns+1), sizeof (long));
    num = calloc ((ncolumns+1), sizeof (long));
    Mat = lrs_alloc_mp_matrix (q+1,ncolumns+1);
    Den = lrs_alloc_mp_vector (ncolumns+1);

    for (i=0; i<=ncolumns; i++)
       {
         den[i]=ONE;
         num[i]=ZERO;
         itomp(ONE,Den[i]);
       }

/* allocate and init structure for static problem data */

    Q = lrs_alloc_dat ("LRS globals");  
    if (Q == NULL)
             return 1;
/* here various flags in lrs_dat Q can be set */

    Q->lponly=TRUE;
    Q->nonnegative=TRUE;  /* non-negative variables */


/*
    Q->debug=TRUE;
    Q->verbose=TRUE;    
*/

/* allocate and initialize lrs_dic                     */

    Q->m=q; Q->n=ncolumns+1;

    P = lrs_alloc_dic (Q);            
    if (P == NULL)
       return 1;


    for(j=1;j<=q;j++) /* set constraints */
     {
      itomp(t[j],temp);
      Binom(l,s[j],Mat[j][0]);
      mulint(temp,Mat[j][0],Mat[j][0]);
     }

    /* the first partition */
    current_column=1;
    for(i=1;i<q;i++)
      a[i]=0;
    a[q]=l-s[0];

    do {

      /* add the current_column's coefficients */

      for(j=1;j<=q;j++)
        {
        Binom(a[j],s[j],Mat[j][current_column]);
        changesign(Mat[j][current_column]);
        }
      /* find next partition */

      i=q+1; /* find the largest i with a[i]>0 */
      while(a[--i]==0) ;
    
      if ( i>1 ) { /* there are more partitions */

	current_column++;
	a[i-1]++; a[q]=a[i]-1;
	if(i<q) a[i]=0; /* beware that i might be equal q */

      }     

    } while ( i>1 );
   

    if(current_column!=ncolumns) { /* something is wrong */
      printf("current_column=%ld ncolumns=%ld l=%ld\n",
	     current_column,ncolumns,l);
      exit(1);
    }

/* set up lrs_dic row by row */
   for (j=1;j<=q;j++)
       lrs_set_row_mp(P,Q,j,Mat[j],Den,GE);
    
/* set up objective function */

    num[0]=0;
    for(j=1;j<=ncolumns;j++)
       num[j]=l;

    lrs_set_obj(P,Q,num,den,MAXIMIZE);

/* now we invoke an lp solver and get its output */

    if(!lrs_solve_lp(P,Q)) {
      printf("Could not solve LP!\n");
      printA(P,Q);
      exit(1);
    }

     printf("\nl=%ld,",l);
     prat (" LP solution=", Q->objnum, Q->objden);
     fflush(stdout);

    if (first_lp) 
      {
      copy(best_bound_num,Q->objnum);
      copy(best_bound_den,Q->objden);
      first_lp=0;
      } 
    else
      if( comprod(Q->objnum,best_bound_den,Q->objden,best_bound_num) == -1)
      {
        copy(best_bound_num,Q->objnum);
        copy(best_bound_den,Q->objden);
        prat ("\n*New best bound ", best_bound_num, best_bound_den);
      }

/* deallocate space   */

    lrs_free_dic (P,Q);          
    lrs_free_dat (Q);    

    lrs_clear_mp_matrix(Mat,q+1,d+1);
    lrs_clear_mp_vector(Den,d+1);

    free(den);
    free(num);

    l++; 
    j=l*t[0];
    itomp (j,lt);
    printf("lt= %ld",j);

  }

  while ( comprod ( lt,best_bound_den,best_bound_num,mpone) != 1 );


  printf("\nThe value of\n $\\lim \\hat r(");
  for(j=1;j<=q;j++)
    printf("K_{%ld,%ldn},", s[j], t[j]);
  for(j=q+1;j<=r;j++)
    printf("K_{%ld,M},",s[j]);

  prat (")/n$\nis ", best_bound_num, best_bound_den);
  lrs_close("lramsey:");
  printf("\n");
  return 0;

}

long
binom(long n, long m)
{

  long i, en=1, den=1;

  if (n<0 || n<m)
    return 0;

  if (2*m>n)
    m=n-m;

  for(i=1;i<=m;i++) {
    en*=n-i+1;
    den*=i;
  }

  return en/den;

}

void Binom(long n, long m, lrs_mp bin)
{
  lrs_mp en, den;
  long i;


  if (n<0 || n<m)
    {
    itomp(ZERO,bin);
    return;
    }

  lrs_alloc_mp (en); lrs_alloc_mp (den);

  itomp(ONE,en);
  itomp(ONE,den);

  if (2*m>n)
    m=n-m;

  for(i=1;i<=m;i++) {
    
    itomp(n-i+1,bin);
    mulint(bin,en,en);
    itomp(i,bin);
    mulint(bin,den,den);
    reduce(en,den);
   }
  divint(en,den,bin);
 
  lrs_clear_mp(en); lrs_clear_mp(den);
  return;
}


\end{verbatim}

\end{document}